\documentclass[a4paper,11pt]{amsart}
\usepackage{amsthm,amssymb,mathtools,amsfonts,latexsym,rawfonts,amsmath, mathrsfs, lscape}
\usepackage{stmaryrd,tikz,pgfplots}
\usepackage{marginnote}  
\usepackage[backref=page]{hyperref}
\usepackage[top=1.5in, bottom=0.9in, left=0.9in, right=0.9in]{geometry}

\usepackage{enumitem}

\newlist{legal}{enumerate}{10}
\setlist[legal]{label*=\arabic*.}
\usepackage[textsize=tiny]{todonotes}

\newcommand{\Marginpar}[1]{\marginpar{\tiny{#1}}}
\newcommand{\Note}[1]{{\par\noindent\hrulefill\par\tiny{#1}\par\noindent\hrulefill\par}}
\newcommand{\Detail}[1]{{#1}}
\renewcommand{\Marginpar}[1]{}
\renewcommand{\Note}[1]{}
\renewcommand{\Detail}[1]{}
\renewcommand*{\backref}[1]{}
\renewcommand*{\backrefalt}[4]{%
    \ifcase #1 (Not cited.)%
    \or        (Cited on page~#2.)%
    \else      (Cited on pages~#2.)%
    \fi}

\usepackage{hyperref}
\usepackage[all]{xy}
\usepackage{amsmath,amsfonts}	
\usepackage{amsthm,latexsym,amssymb}
\usepackage{tensor}

\usepackage{array, tabularx}

\usepackage{booktabs} 
\usepackage{siunitx}

\usepackage{tikz-cd}

\newtheorem{thm}{Theorem}
\newtheorem*{thm*}{Theorem}
\newtheorem{prop}[thm]{Proposition}
\newtheorem*{prop*}{Proposition}

\newtheorem*{lem*}{Lemma}

\newtheorem*{cor*}{Corollary}

\newtheorem* {Question*}{Question}

\theoremstyle{definition}
\newtheorem{defn}[thm]{Definition}
\newtheorem*{defn*}{Definition}

\renewcommand{\[}{\begin{equation*}}
\renewcommand{\]}{\end{equation*}}

\DeclareMathOperator\vol{vol}

\def \R {\mathbb R}

\begin{document}
\parskip 1mm

\title[Pluricanonical metrics and circle-invariant $G_2$
structures]{Pluricanonical locally conformally almost K\"ahler metrics,\\
circle-invariant $G_2$ structures, and Laplacian flows.}

\author{Mehdi Lejmi}
\address{Department of Mathematics, Bronx Community College of CUNY, Bronx, NY 10453, USA.}
\email{mehdi.lejmi@bcc.cuny.edu}

\thanks{ }



\maketitle
\begin{abstract}
We study explicit Hodge evolutions of pluricanonical locally conformally
almost K\"ahler metrics and their relation to locally conformal parallel
$G_2$ structures. A compact locally conformal parallel $7$-manifold with
a parallel non-zero Lee form, carrying a structure-preserving free circle action
of constant orbit length orthogonal to the Lee field, has a $6$-dimensional
quotient which is pluricanonical with a parallel Lee form and a non-integrable
almost complex structure. The fundamental forms upstairs and downstairs
are Laplacian eigenforms with an explicit eigenvalue ratio. Conversely,
for a fixed normalized $SU(3)$ completion, we characterize the curvature
of a circle-invariant locally conformal parallel lift. The lift equations
force a pluricanonical base to have parallel Lee form, even locally, and
exclude strict integrable pluricanonical bases. Independently, the
pluricanonical eigenform identity gives a shrinking solution of He's full
Hodge equation and an immortal $d^*d$ solution contracting the Lee plane.

\end{abstract}

\section{Introduction}
A locally conformally symplectic (LCS) structure on a manifold $M$ of real
dimension $2n\geq4$ is a pair $(F,\theta)$, where $F$ is a nondegenerate
$2$-form and $\theta$ is a closed $1$-form such that $$dF=\theta\wedge F.$$
The form $\theta$ is called the Lee form. If $F$ is compatible with an
almost complex structure $J$, the Riemannian metric $g=F(\cdot,J\cdot)$
is locally conformally almost K\"ahler (LCAK). When $J$ is integrable,
the metric is locally conformally K\"ahler (LCK). An LCK metric with
parallel Lee form is called Vaisman. 

Pluricanonical LCAK metrics were introduced in \cite{addison2026pluricanonicallocallyconformallykahler} as an extension
of the pluricanonical condition considered by Kokarev in LCK geometry
\cite{MR2520354}. The condition requires $D\theta$ to be $J$-anti-invariant
and the image of the Nijenhuis tensor to be orthogonal to the Lee field
$T=\theta^\sharp$ and the anti-Lee field $JT$ (here $D$ is the Levi-Civita connection). On compact complex
manifolds, pluricanonical LCK metrics are Vaisman \cite{MR3674175}.
In the non-integrable case the Lee form need not be parallel, even on
compact manifolds. We call a pluricanonical LCAK structure \emph{strict} if
$\theta\not\equiv0$, and \emph{almost Vaisman} if $D\theta=0$.

One of the identities established for pluricanonical LCAK metrics in \cite{addison2026pluricanonicallocallyconformallykahler} is
\begin{equation}\label{eq:eigenform}
 \Delta^gF=(n-1)\|\theta\|^2F,
\end{equation}
where $\Delta^g=d\delta^g+\delta^gd$ is the Hodge Laplacian
and $\|\theta\|^2=g(\theta,\theta)$ is constant (here $\delta^g$ is the codifferential). This places the
fundamental $2$-form in a distinguished class for geometric evolutions.
It also suggests asking whether pluricanonical structures can arise
naturally as quotients of seven-dimensional geometries with analogous
differential identities. We address both questions, keeping track of
the metric used by the Laplacian and of the distinction between
integrable and non-integrable almost complex structures.
Our geometric starting point is a positive $3$-form $\varphi$ with
\begin{equation}\label{eq:intro-lcp}
 d\varphi=\vartheta\wedge\varphi,\qquad
 d(*_\varphi\varphi)=\frac43\vartheta\wedge *_\varphi\varphi,\qquad
 d\vartheta=0,
\end{equation}
where $\vartheta$ is the Lee
form of $\varphi$ and $ *_\varphi$ is the Hodge operator induced by $\varphi.$
These are the locally conformal parallel equations in our Lee
normalization. Ivanov--Parton--Piccinni \cite{gtwo-IPP2006} describe compact
structures in a Gauduchon gauge and their nearly K\"ahler product
covers. A structure-preserving circle action supplies an almost
Hermitian quotient. The question relevant here is whether that quotient
satisfies the full pluricanonical condition, including the Nijenhuis
orthogonality requirement. The following theorem gives an affirmative answer.


\begin{thm*}[Theorem~\ref{gtwo-recognition}]
Let $M^7$ be compact and connected, with a locally conformal parallel
$G_2$ structure $\varphi$ whose nonzero Lee form $\vartheta$ is parallel.
Suppose that a free circle action preserves $\varphi$ and its generator
$X$ has constant length $c>0$ and satisfies $\vartheta(X)=0$.
Then $B=M/S^1$, with its quotient metric, carries a strict almost
Vaisman structure with a non-integrable almost complex structure, given by
\[
 \pi^*F=\iota_{X/c}\varphi,\qquad \pi^*\theta=\vartheta.
\]
Moreover,
\[
 \Delta^{g_B}F=2\|\vartheta\|^2F,\qquad
 \Delta^{g_\varphi}\varphi=\frac43\|\vartheta\|^2\varphi.
\]
The universal cover of $M$ is isometric to a product $\R\times N$, where
$N$ is homothetic to $S^3\times S^3$ with its canonical nearly K\"ahler
structure.
\end{thm*}

The converse problem has additional differential and topological
constraints.
First, we fix a normalized $\mathrm{SU}(3)$
completion $(F,\psi_+,\psi_-)$ of a $6$-dimensional LCAK structure.
Proposition~\ref{gtwo-completion-criterion} identifies the circle-lift equations
as
\begin{equation}\label{eq:intro-system}
 \kappa\wedge F+d\psi_+=\theta\wedge\psi_+,\qquad
 d\psi_-=\frac43\theta\wedge\psi_-,\qquad
 \kappa\wedge\psi_-=\frac13\theta\wedge F^2.
\end{equation}
Here $\kappa$ is the candidate curvature of the unit vertical $1$-form.
For a circle parameter of period $2\pi$ and fixed radius $c>0$, a global
circle lift exists precisely when $\kappa$ is closed and $[\frac{\kappa}{2\pi c}]$ is integral. 
In particular, the topological triviality of the canonical line bundle
alone does not solve the lift problem. The circle-bundle differentiation formulas underlying
\eqref{eq:intro-system} belong to the established $SU(3)$--$G_2$
framework of Chiossi--Salamon \cite{MR1922042}.
The lift equations
impose the following additional rigidity on a pluricanonical LCAK base.

\begin{thm*}[Theorem~\ref{gtwo-local-parallel-obstruction}]
If a pluricanonical LCAK $6$-manifold has a normalized $\mathrm{SU}(3)$
completion and a closed curvature form $\kappa$ satisfying the circle-lift
equations, then
\[
 D\theta=0,\qquad\iota_T\kappa=0,
\]
where $T=\theta^\sharp.$
\end{thm*}
Its proof uses the constant length and coclosedness of the
lifted Lee form, conformal Ricci flatness, and a pointwise Bochner
identity (see \cite{MR2252869} for instance).
The curvature-norm identity in Proposition~\ref{gtwo-norm-obstruction}
also shows that an integrable pluricanonical base satisfying these
equations with $d\kappa=0$ has $\theta=\kappa=0$. Then, we give an explicit compact
lift from $S^2\times S^3\times S^1$, obtained by reducing the product
of the canonical nearly K\"ahler $S^3\times S^3$ with a circle, to illustrate the pluricanonical quotient
and the integral curvature condition.

Moreover, the eigenform identity~\eqref{eq:eigenform}  leads to explicit solutions of the two
flows introduced by He \cite{flow-He2021}. On a connected manifold, the full Hodge-Laplacian flow shrinks
homothetically, whereas the $d^*d$ flow fixes the horizontal part of the
fundamental form and contracts the Lee plane.

\begin{thm*}[Theorems~\ref{flow-homothety} and~\ref{flow-immortal-ddstar}]
Let $(g_0,J,F_0,\theta)$ be a strict pluricanonical LCAK structure in
real dimension $2n\geq4$. Set
\[
 r^2=\|\theta\|_{g_0}^2,\qquad \lambda=(n-1)r^2,
 \qquad F_V=r^{-2}\theta\wedge J\theta,\qquad F_H=F_0-F_V.
\]
With $J$ fixed and $g(t)=F(t)(\cdot,J\cdot)$, the two families
\[
 \begin{array}{lll}
 F(t)=(1-\lambda t)F_0
 &\text{solve }\partial_tF=-\Delta_{g(t)}F,
 &t<\lambda^{-1},\\[3pt]
 F(t)=F_H+(1+\lambda t)^{-1}F_V
 &\text{solve }\partial_tF=-\delta^{g(t)}dF,
 &t>-\lambda^{-1},
 \end{array}
\]
with initial form $F_0$, and remain pluricanonical. Their Lee forms are
$\theta$ and $(1+\lambda t)^{-1}\theta$, respectively. On a compact
manifold, the second family is the unique smooth forward solution in
He's fixed-$J$ tamed class with initial form $F_0$, on every common
interval of existence.
\end{thm*}

Thus the first solution degenerates at finite positive time, while the
second is immortal and converges to the closed rank-$(2n-2)$ form $F_H$.
Similar deformations of Vaisman metrics were considered for instance in~\cite{MR3342236} (see also the references therein).

The paper is organized as follows. In Section~\ref{sec:preliminaries}, we fix some notations and conventions. Then, in Section~\ref{gtwo-circle}, we recognize the quotient of a compact locally conformal parallel
$G_2$ structure, carrying a structure-preserving free circle action, as a pluricanonical structure. Moreover, we study the circle-lift equations on a pluricanonical manifold to a locally conformal parallel $G_2$ structure. Finally, in section~\ref{flow-he}, we use the eigenform property of pluricanonical LCAK metrics to solve He's Laplacian flows.

{\bf{Acknowledgement}}: The author is very grateful to Tedi Dr\u{a}ghici for the discussion about Donaldson's idea~\cite{MR3702382} to connect hypersymplectic structures on $4$-manifolds and $G_2$ structures, and for suggesting the direct first-kind LCS construction used in the proof of Theorem~\ref{gtwo-recognition}.

{\bf{Declaration of AI assistance}}: The author acknowledges the use of ChatGPT 6.1 Sol (OpenAI) for checking routine calculations and improving manuscript exposition. All mathematical ideas, proof strategies, and results are the sole work of the author, who independently verified every step of the proofs.


\section{Preliminaries}\label{sec:preliminaries}

Let $(M^{2n},F,\theta)$ be a locally conformally symplectic (LCS for short) manifold of real dimension $2n$, where $F$ is a non-degenerate $2$-form and $\theta$ is a closed $1$-form such that $dF=\theta\wedge F.$ The $1$-form $\theta$ is called the Lee form of the LCS structure. We suppose that $F$ is compatible with an almost complex structure $J$, that is $g:=F(\cdot,J\cdot)$ is a Riemannian metric. The induced metric $g$ is called a locally conformally almost K\"ahler (LCAK for short) metric. We denote by $D$ the Levi-Civita connection of $g.$ 
We use the convention $(J\alpha)(X)=-\alpha(JX)$ for a $1$-form $\alpha$, and $\operatorname{vol}_g=\frac{F^n}{n!}$. For a covariant $2$-tensor $\phi$, its $J$-invariant and $J$-anti-invariant
parts are $\phi^{J,\pm}(X,Y)=\frac12(\phi(X,Y)\pm \phi(JX,JY))$. Our Nijenhuis tensor convention is
\[
4N(X,Y)=[JX,JY]-[X,Y]-J[JX,Y]-J[X,JY].
\]
The Hodge Laplacian is defined as $\Delta_g=d\delta_g+\delta_gd$, where $\delta^g$ is the codifferential, defined as the adjoint of the exterior derivative $d$ with respect to the metric $g$.
\begin{defn}\label{def:pc}~\cite{addison2026pluricanonicallocallyconformallykahler}
An LCAK metric is pluricanonical if $D\theta$ is $J$-anti-invariant and the image of the Nijenhuis tensor is $g$-orthogonal to $T=\theta^\sharp$ and $JT$. A pluricanonical LCAK structure is strict if its Lee form is not identically zero. 
\end{defn}

For such a structure, in addition to~\eqref{eq:eigenform}, we also have
\begin{equation}\label{intro-pc-identities}
\delta^g\theta=0,\quad d\|\theta\|^2=0,\quad
dJ\theta=-\|\theta\|^2F+\theta\wedge J\theta,\quad
\delta^gF=-(n-1)J\theta,\quad \Delta^g J \theta=(n-1)\|\theta\|^2J\theta.
\end{equation}
Moreover, we recall that for a given pluricanonical LCAK structure, the tangent bundle splits into
\begin{equation}\label{splitting}
TM=\mathrm{span}\{T,JT\}\oplus H,
\end{equation}
where $H=\mathrm{span}\{T,JT\}^{\perp_g}.$

\section{Locally conformal parallel \texorpdfstring{$G_2$}{G2} structures}\label{gtwo-circle}

\subsection{Recognition of pluricanonical quotients}

Let $\varphi$ be a positive $3$-form on a $7$-manifold. Namely, if we choose an orientation and any positive volume form $\nu$, then the symmetric bilinear form $B$ defined by
$$B_\varphi(X,Y)\nu=\iota_X\varphi\wedge \iota_Y\varphi\wedge\varphi$$
is positive definite, where $X,Y$ are vector fields and $\iota_X$ is the contraction by $X.$
The metric $g_\varphi$ induced by $\varphi$ is characterized then by
$$g_\varphi(X,Y)\nu_{g_\varphi}=\frac{1}{6}\iota_X\varphi\wedge \iota_Y\varphi\wedge\varphi.$$
The dual $4$-form $\psi$ is defined as $$\psi=*_{\varphi}\varphi.$$
The $3$-form $\varphi$ is called torsion-free if $D^{g_\varphi}\varphi=0$ which is equivalent to $d\varphi=d\psi=0$ (here $D^{g_\varphi}$ is the Levi-Civita connection of the metric $g_\varphi$). Moreover,
 $\varphi$ is called locally conformally parallel if the following equations are satisfied
 \begin{equation}\label{gtwo-lcp}
 d\varphi=\vartheta\wedge\varphi,\qquad
 d\psi=\frac43\vartheta\wedge\psi,\qquad d\vartheta=0,
\end{equation}
where the closed $1$-form $\vartheta$ is the Lee form. For a compact locally conformal parallel $G_2$ structure, the Gauduchon representative has a parallel Lee form \cite{gtwo-IPP2006}. Locally, if $\vartheta=dh$, the form $e^{-h}\varphi$ is torsion-free and its metric is $e^{-2h/3}g_\varphi$. It turns out that a circle quotient of a locally conformal parallel structure yields a pluricanonical base under natural metric assumptions.

\begin{thm}\label{gtwo-recognition}
Let $(M,\varphi)$ be a $7$-dimensional compact and connected $G_2$ manifold satisfying \eqref{gtwo-lcp} with $\vartheta\ne0$ and $D^{g_\varphi}\vartheta=0$. Suppose a free circle action preserves $\varphi$, its generator $X$ has constant length $c>0$, and $\vartheta(X)=0$. Set $V=X/c$. Then the quotient metric and the form
\[
 \pi^*F=\iota_V\varphi
\]
define an almost Vaisman structure $(g_B,J_B,F,\theta)$ on $B=M/S^1$, where $\pi^*\theta=\vartheta$, and $\pi:M\rightarrow B$ is the quotient map. Furthermore, we have
\begin{equation}\label{gtwo-quotient-eigenvalues}
 \Delta^{g_B}F=2\|\vartheta\|^2F,\qquad
 \Delta^{g_\varphi}\varphi=\frac{4\|\vartheta\|^2}{3}\varphi.
\end{equation}
The induced almost complex structure $J_B$ is non-integrable. Moreover, the universal cover of $M$ has a nearly K\"ahler factor homothetic to the canonical $S^3\times S^3$.
\end{thm}
\begin{proof}


Since the circle action preserves $\varphi$, it preserves the metric
$g_\varphi$. Let $\pi:M\to B=M/S^1$ be the quotient map. The horizontal, invariant
$2$-form $\iota_V\varphi$ descends to $F$ on $B$.
Pointwise $G_2$ algebra gives an almost complex structure $J_B$
compatible with the quotient metric $g_B$, with
$F(Y,Z)=g_B(J_BY,Z)$. The closed form $\vartheta$ descends also on $B$ because
$\vartheta(V)=0$. We write $\pi^*\theta=\vartheta$. In the following, we suppress pullbacks when
no confusion can arise.

Locally, $\vartheta=dh$ and $e^{-h}\varphi$ is parallel for
$e^{-2h/3}g_\varphi$. The conformal connection formula therefore gives the global
tensor identity
\begin{equation}\label{covariant-derv-G2}
D^{g_\varphi}_Y\varphi=\frac13\left(Y^\flat\wedge\iota_{\vartheta^\sharp}\varphi
                     -\vartheta\wedge\iota_Y\varphi\right).
\end{equation}
Since $D^{g_\varphi}\vartheta^\sharp=0$, this implies $D^{g_\varphi}_{\vartheta^\sharp}\varphi=\mathcal{L}_{\vartheta^\sharp}\varphi=0$. Moreover,
$\mathcal{L}_V\vartheta=0$ and $\mathcal{L}_Vg_\varphi=0$, so $[V,\vartheta^\sharp]=0$. Thus $\vartheta^\sharp$ is horizontal
and projects to $T=\theta^\sharp$. For basic horizontal vector fields $Y,Z$,
the horizontal connection formula gives
\[
g_B(D^{g_B}_{\pi_*Y}T,\pi_*Z)=g(D^{g_\varphi}_Y\vartheta^\sharp,Z)=0.
\]
Consequently $D^{g_B}\theta=0$.

Cartan's formula gives
\[
dF=-\iota_Vd\varphi=\theta\wedge F.
\]
Let $S=\iota_{\vartheta^\sharp}\varphi$. Then
\[
dS=-\|\vartheta\|^2\varphi+\vartheta\wedge S,\qquad
\pi^*(J_B\theta)=\iota_{\vartheta^\sharp}\iota_V\varphi=-\iota_VS.
\]
Since $\mathcal{L}_VS=0$, differentiation yields
\[
dJ_B\theta=-\|\vartheta\|^2F+\theta\wedge J_B\theta.
\]
Also, $\mathcal{L}_{\vartheta^\sharp}F=0$ upstairs, and hence $\mathcal{L}_TF=0$ downstairs.
The form $dJ_B\theta$ is $J_B$-invariant. Since
$d\theta=0$, it follows that
\[
\begin{aligned}
4\theta(N(Y,Z))
 &=-dJ_B\theta(J_BY,Z)-dJ_B\theta(Y,J_BZ)=0,\\
4J_B\theta(N(Y,Z))
 &=dJ_B\theta(Y,Z)-dJ_B\theta(J_BY,J_BZ)=0.
\end{aligned}
\]
Together with $D^{g_B}\theta=0$, these identities show that the quotient
is strict almost Vaisman. Then, the first equation in~\eqref{gtwo-quotient-eigenvalues} follows from Equation~\eqref{eq:eigenform} in dimension $2n=6.$

To get the second equation in~\eqref{gtwo-quotient-eigenvalues}, we first trace Equation~\eqref{covariant-derv-G2}. We get
$\delta^{g_\varphi}\varphi=-\frac43S$. Since $D^{g_\varphi}_{\vartheta^\sharp}\varphi=0$, we obtain
\[
\begin{aligned}
\Delta^{g_\varphi}\varphi
 &=-\frac43dS+\delta^{g_\varphi}(\vartheta\wedge\varphi)\\
 &=\frac43\|\vartheta\|^2\varphi-\frac43\vartheta\wedge S
                         +\frac43\vartheta\wedge S
 =\frac43\|\vartheta\|^2\varphi.
\end{aligned}
\]

To prove the non-integrability, we use the natural normalized $SU(3)$ completion
on $B$ defined by
\[
\psi_+=\varphi-V^\flat\wedge F,\qquad
\psi_-=-\iota_V(*_{g_\varphi}\varphi).
\]
These forms are horizontal and invariant, so they descend on $B$. Pointwise $G_2$ algebra gives
\[
*_{g_\varphi}\varphi=\frac12F^2-V^\flat\wedge\psi_-,
\qquad \Psi=\psi_++\sqrt{-1}\psi_-\in\Lambda^{3,0}(B).
\]
Comparing the vertical and horizontal parts of
$d(*_{g_\varphi}\varphi)=\frac43\vartheta\wedge *_{g_\varphi}\varphi$ gives
\[
d\psi_-=\frac43\theta\wedge\psi_-,
\qquad
dV^\flat\wedge\psi_-=\frac13\theta\wedge F^2.
\]
The second identity implies
$(dV^\flat)^{J_B,-}=\frac13\iota_T\psi_-$.
Indeed, wedging with $\psi_-$ is an isomorphism on real
$J_B$-anti-invariant $2$-forms, while $J_B$-invariant $2$-forms wedge
trivially with $\psi_\pm$. The standard $SU(3)$ identities
\[
(\iota_T\psi_-)\wedge\psi_-=\theta\wedge F^2,\qquad
(\iota_T\psi_-)\wedge\psi_+=-2*_{g_B}\theta
\]
therefore give
\begin{equation}\label{eq:psi+}
d V^\flat\wedge\psi_+=-\frac23*_{g_B}\theta.
\end{equation}

Suppose that $J_B$ is integrable. The complex $4$-form
\[
Q=d\Psi-\frac43\theta\wedge\Psi
\]
is then of type $(3,1)$. The equation for $d\psi_-$ makes $Q$ real,
so $Q=0$. In particular,
\begin{equation}\label{eq:dpsi+}
d\psi_+=\frac43\theta\wedge\psi_+.
\end{equation}
Since
$D^{g_B}\theta=0$, we obtain using Equations~\eqref{eq:psi+} and \eqref{eq:dpsi+}
\[
\begin{aligned}
0
 &=-\frac23d(*_{g_B}\theta)
   =d(dV^\flat\wedge\psi_+)\\
 &=dV^\flat\wedge d\psi_+
   =\frac43\theta\wedge dV^\flat\wedge\psi_+
   =-\frac89\|\vartheta\|^2\vol_{g_B},
\end{aligned}
\]
a contradiction. Hence $J_B$ is non-integrable.

It remains only to identify the global nearly K\"ahler factor. By
\cite[Proposition 3.3 and Corollary 3.4]{gtwo-IPP2006}, the universal cover of $M$
is a Riemannian product $\widetilde M=\mathbb R\times N$, with
\begin{equation}\label{eq:product}
g_\varphi=dt^2+g_N,\qquad \vartheta=-\|\vartheta\|\,dt,\qquad
\varphi=dt\wedge\omega_N+\Psi_{N,+},
\end{equation}
where $N$ is a compact, simply connected, and strictly nearly K\"ahler manifold.
Let $E=\partial_t$. Then, the induced forms have the normalization
\begin{equation}\label{eq:nk}
\begin{gathered}
\omega_N=\iota_E\varphi,\qquad
d\omega_N=\|\vartheta\|\Psi_{N,+},\qquad d\Psi_{N,+}=0,\\
d\Psi_{N,-}=-\frac{2\|\vartheta\|}{3}\omega_N^2,\qquad
\Psi_{N,-}=*_{g_N}\Psi_{N,+},\qquad
\omega_N\wedge\Psi_{N,+}=0.
\end{gathered}
\end{equation}
The lift of $V$ is tangent to $N$ because $\vartheta(V)=0$.
Since $\vartheta^\sharp=-\|\vartheta\|\,\partial_t$ and $[\vartheta^\sharp,V]=0$, this lifted field is independent
of $t$ and restricts to a unit Killing field on $N$.
After a constant rescaling, \cite[Theorem 1.1]{gtwo-MNS2005} applies. As $N$ is
simply connected, it is homothetic to $S^3\times S^3$ with its canonical
nearly K\"ahler structure.
\end{proof}

\subsection{The circle-lift equations}

Let $(B,g,J,F,\theta)$ be a $6$-dimensional LCAK manifold. Suppose that a normalized $SU(3)$ completion $(F,\psi_+,\psi_-)$ is fixed, with
\begin{equation}\label{gtwo-su3}
 F\wedge\psi_\pm=0,\qquad
 \psi_+\wedge\psi_-=\frac23F^3,\qquad *_g\psi_+=\psi_-.
\end{equation}
In particular, $\Psi=\psi_++i\psi_-$ is  a nowhere vanishing normalized complex volume $3$-form. Such an $SU(3)$ completion exists when the (smooth) canonical line bundle $\Lambda^{3,0}(B)$ is trivial which is equivalent to $c_1(B,J)=0\in H^2(B,\mathbb{Z}).$ Let $\pi:P\longrightarrow B$ be some principal circle bundle, whose principal generator $X$ has period $2\pi$. Fix a connection $A$ and a constant $c>0$, and set $\eta=c A$.
The connection metric is $g_P=\pi^*g+\eta^2$, and $X$ has length $c$ (the connection $A$ is normalized by $A(X)=1$).
Define the following forms (pullbacks are suppressed below):
\begin{equation}\label{gtwo-circle-forms}
 \varphi=\eta\wedge F+\psi_+,\qquad
 \psi=\ast_{g_P}\varphi=\frac12F^2-\eta\wedge\psi_-.
\end{equation}
These formulas define a positive $G_2$ structure with metric $g_P$ (this can be checked pointwise).
\begin{prop}\label{gtwo-completion-criterion}
For the fixed completion \eqref{gtwo-su3}, the forms \eqref{gtwo-circle-forms} satisfy \eqref{gtwo-lcp}, with $\vartheta=\pi^*\theta$, if and only if
\begin{align}
 \kappa\wedge F+d\psi_+&=\theta\wedge\psi_+,\label{gtwo-lift1}\\
 d\psi_-&=\frac43\theta\wedge\psi_-,\label{gtwo-lift2}\\
 \kappa\wedge\psi_-&=\frac13\theta\wedge F^2.\label{gtwo-lift3}
\end{align}
Equivalently, writing $T=\theta^\sharp$ and
\[
 \mathcal D=d\psi_+-\frac43\theta\wedge\psi_+,
\]
one must have
\begin{equation}\label{gtwo-curvature-decomposition}
 d\psi_-=\frac43\theta\wedge\psi_-,\qquad
 \mathcal D\in\Lambda^{2,2}_0(B),\qquad
 \kappa=\frac13\iota_T\psi_-+\kappa^{J,+}_0,\quad \kappa^{J,+}_0=*_g\mathcal D,
\end{equation}
where $\Lambda^{2,2}_0(B)=\{D\in \Lambda^{2,2}(B)\,|\, D\wedge F=0 \}$ and $\kappa^{J,+}_0$ is the trace-free $J$-invariant part of $\kappa.$ The candidate curvature $\kappa$ is unique. A global circle lift with principal radius $c$ exists precisely when this curvature is closed and
\begin{equation}\label{gtwo-integrality}
 \left[\frac{\kappa}{2\pi c}\right]\in
 \operatorname{im}\bigl(H^2(B,\mathbb Z)\longrightarrow H^2(B,\mathbb R)\bigr).
\end{equation}
\end{prop}
\begin{proof}
Differentiating $\psi$ in \eqref{gtwo-circle-forms}, we obtain
$$d\psi=\theta\wedge F^2-\kappa\wedge\psi_-+\eta\wedge d\psi_-.$$
Comparing the horizontal and vertical components to
$$\frac{4}{3}\theta\wedge\psi=\frac{2}{3}\theta\wedge F^2+\frac{4}{3}\eta\wedge\theta\wedge\psi_-,$$
we obtain \eqref{gtwo-lift2} and \eqref{gtwo-lift3}. In a similar way, to obtain \eqref{gtwo-lift1}, we differentiate $\varphi$ in \eqref{gtwo-circle-forms}.
 
Moreover, the map $L_F:\Lambda^2(B)\longrightarrow\Lambda^4(B)$, $\zeta\longmapsto\zeta\wedge F$, is an isomorphism in dimension six, so \eqref{gtwo-lift1} determines $\kappa$ uniquely.

Furthermore, differentiating $F\wedge \psi_+=0$ gives that $F\wedge d\psi_+=0$, and hence $\kappa\wedge F^2=0$. Decompose
$$\kappa=\kappa^{J,-}+\kappa^{J,+}_0.$$
 In an orthonormal $SU(3)$ coframe one has
\begin{equation}\label{gtwo-algebra}
 (\iota_T\psi_-)\wedge\psi_-=\theta\wedge F^2,\quad
 (\iota_T\psi_-)\wedge F=-\theta\wedge\psi_+,\quad
 (\iota_T\psi_-)\wedge\psi_+=-2*_g\theta.
\end{equation}
Wedging with $\psi_-$ is an isomorphism from the real $J$-anti-invariant $2$-forms to $\Lambda^5(B)$, whereas $\kappa^{J,+}_0\wedge\psi_\pm=0.$
Thus \eqref{gtwo-lift3} gives $3\kappa^{J,-}=\iota_T\psi_-$. Substitution into \eqref{gtwo-lift1} gives $\mathcal D=-\kappa^{J,+}_0\wedge F=*_g\kappa^{J,+}_0$, proving \eqref{gtwo-curvature-decomposition}. The converse follows from the same computations.

The condition \eqref{gtwo-integrality} is the Chern class condition given that the curvature of $A$ is $\frac{\kappa}{c}$.

\end{proof}

\begin{prop}\label{gtwo-norm-obstruction}
If the system \eqref{gtwo-lift1}--\eqref{gtwo-lift3} has a closed curvature, then
\begin{equation}\label{gtwo-lee-curvature-norm}
 \delta^g\theta=-\frac43\|\theta\|^2+\frac32\|\kappa^{J,+}_0\|^2.
\end{equation}
On a pluricanonical base this implies
\[
 \|\kappa^{J,+}_0\|^2=\frac89\|\theta\|^2,\qquad
 \|\kappa\|^2=\frac{10}{9}\|\theta\|^2.
\]
If the base is integrable and pluricanonical, then $\theta=0$ and $\kappa=0$. More generally, the same conclusion holds for a compact LCK base.
\end{prop}
\begin{proof}
By \eqref{gtwo-algebra}, $\kappa\wedge\psi_+=-\frac{2}{3}*_g\theta$. Recall that
$d\psi_+=\frac{4}{3}\theta\wedge\psi_++\ast_g \kappa^{J,+}_0.$
Since $d\kappa=0$, we obtain
\begin{eqnarray*}
 d(\kappa\wedge\psi_+)&=&\kappa\wedge d\psi_+\\
&=&\frac{4}{3}\kappa\wedge \theta\wedge\psi_++\kappa\wedge\ast_g \kappa^{J,+}_0\\
&=&-\frac89\theta\wedge*_g\theta+\kappa\wedge\ast_g \kappa^{J,+}_0\\ 
&=&\left(-\frac89\|\theta\|^2+\|\kappa^{J,+}_0\|^2\right)\operatorname{vol}_g.
\end{eqnarray*}
But $d(\kappa\wedge\psi_+)=-\frac{2}{3}d*_g\theta=\frac{2}{3}\delta^g\theta\,\operatorname{vol}_g$, so we get \eqref{gtwo-lee-curvature-norm}.

When $J$ is integrable, $d\Psi-\frac43\theta\wedge\Psi$ has type $(3,1)$. By \eqref{gtwo-curvature-decomposition} it is also the real $(2,2)$ form $- \kappa^{J,+}_0\wedge F$, so it vanishes and $ \kappa^{J,+}_0=0$. Equation \eqref{gtwo-lee-curvature-norm} then gives $\theta=0$ on a pluricanonical base. For a compact integrable base, its integral gives the same conclusion. Finally \eqref{gtwo-curvature-decomposition} gives $\kappa=0$.

\end{proof}

\begin{thm}\label{gtwo-local-parallel-obstruction}
Suppose a pluricanonical LCAK $6$-dimensional manifold  $(B,g,J,F,\theta)$ admits an $SU(3)$  completion and a closed curvature solving \eqref{gtwo-lift1}--\eqref{gtwo-lift3}. Then
\[
 D^g\theta=0,\qquad \iota_T\kappa=0.
\]
\end{thm}
\begin{proof}
Locally write $\kappa=dA$ and take the lift $\eta=ds+A$ on a product with an interval. The metric and positive form are given by \eqref{gtwo-circle-forms}. Write $\Theta=\pi^*\theta$ and $\mathcal T=\Theta^\sharp$. The unit fibres are Killing and geodesic, so $\delta^{g_P}\Theta=\pi^*(\delta^g\theta)=0$. Also $d\Theta=0$ and $\|\Theta\|=r$ is constant by pluricanonicity.

For $dh=\Theta$, the form $e^{-h}\varphi$ is torsion-free. Its metric $\widetilde g=e^{-2h/3}g_P$ is Ricci-flat \cite[Corollaries 2.28 and 4.12]{gtwo-Karigiannis2009}. The conformal Ricci formula in dimension seven therefore gives
\[
 \operatorname{Ric}_{g_P}=-\frac53D^P\Theta-\frac59\Theta\otimes\Theta
 +\frac13 \left(\delta^{g_P}\Theta\right) g_P+\frac59r^2g_P.
\]
Since $\Theta$ is of constant length, we have $(D^P\Theta)(\mathcal T,\mathcal T)=0$, and consequently $\operatorname{Ric}_{g_P}(\mathcal T,\mathcal T)=0$. The pointwise Bochner formula for the harmonic, constant-length form $\Theta$ now yields $\|D^P\Theta\|^2=0$.

For horizontal lifts $Y,Z$ and the unit vertical vector $V$, the connection metric formulas are
\[
 (D^P\Theta)(Y,Z)=(D^B\theta)(Y,Z),\qquad
 (D^P\Theta)(Y,V)=\frac12\kappa(T,Y).
\]
Both conclusions follow.
\end{proof}

We give now an explicit compact solution of the lift system. The quaternionic model and its circle symmetries are described in \cite[Section 2]{gtwo-Dixon2019}. On the other hand, the circle reduction of the cone over the homogeneous nearly K\"ahler space is studied in \cite[Section 3]{gtwo-SaEarpStein2026}. On $S^3\times S^3=SU(2)\times SU(2)$, we write $a=p^{-1}dp$ and $b=q^{-1}dq$ for the imaginary quaternion-valued Maurer--Cartan forms. Thus $da_i=-2a_j\wedge a_k$ and $db_i=-2b_j\wedge b_k$ for cyclic $(i,j,k)$. The canonical nearly K\"ahler structure is
$$
 g_N=\sum_{i=1}^3(a_i\otimes a_i-\frac{1}{2}\left(a_i\otimes b_i+b_i\otimes a_i\right)+b_i\otimes b_i),\qquad
 \omega_N=\frac{\sqrt3}{2}\sum_{i=1}^3a_i\wedge b_i,\qquad
 \Omega=\bigwedge_{i=1}^3(a_i+e^{\frac{2\pi \sqrt{-1}}{3}}b_i),$$$$
 \Psi_{N,+}=-\operatorname{Im}\Omega,\qquad
 \Psi_{N,-}=\operatorname{Re}\Omega.
$$
If we differentiate the above equations, we get \eqref{eq:nk} with $\|\vartheta\|=2$.

We will take a circle quotient of $$P=S^3\times S^3\times(\mathbb R/\ell\mathbb Z),\qquad \ell>0.$$
Fix a unit imaginary quaternion $I$. The quotient map $\pi:P\rightarrow B$ is $$\pi(p,q,t)=(p^{-1}Ip,q,t)\in B=S^2\times S^3\times S^1.$$
On $P$, we consider the $G_2$-structure $\varphi=dt\wedge\omega_N+\Psi_{N,+}$ inducing the metric $g_\varphi=dt^2+g_N.$ Thus, the Lee form $\vartheta=-2dt,$ but the circle action leaves $t$ unchanged so $$\theta_0=-2dt.$$ Moreover, we can find the unit vertical vector and its dual $1$-form. The circle action is $(p,q,t)\longmapsto(e^{Is}p,q,t)$, $s\in\mathbb R/(2\pi\mathbb Z)$. It preserves the nearly K\"ahler structure and has a unit-length generator $V=(Ip,0,0)$. Hence, $$a(V)=p^{-1}Ip=x\in S^2,\qquad b(V)=0,\qquad dt(V)=0.$$
Remark that $x$ is unchanged under the circle action. Moreover, writing the nearly K\"ahler metric as $g_N=|a-\frac{1}{2}b|^2+\frac{3}{4}|b|^2$,
(here $|a|^2=\sum_{i=1}^3a_i\otimes a_i$ etc) we obtain that $g_N(V,V)=1$. Furthermore, for any tangent vector $Y$, we get $g_N(V,Y)=x\cdot\left(a(Y)-\frac{1}{2}b(Y)\right)$ (here $\cdot$ is the Euclidean product). Define
$$A_0=x\cdot a,\qquad \xi=x\cdot b.$$
Then, $$\eta=V^\flat=x\cdot a-\frac{1}{2}x\cdot b=A_0-\frac{1}{2}\xi.$$
Differentiating $x=p^{-1}Ip$ and using that $dp=pa$, we get that $$a=x\,A_0-\frac{1}{2}x\times dx.$$ Then, we can write
$$a-\frac{1}{2}b=x\eta+ \left(-\frac12x\times dx-\frac12(b-x\xi)\right).$$
We remark that the vector-valued $1$-form $-\frac12x\times dx-\frac12(b-x\xi)$ is perpendicular to $x$.
Thus,
$$|a-\frac{1}{2}b|^2=|\eta|^2+\left|-\frac12x\times dx-\frac12(b-x\xi)\right|^2.$$
In particular,
$$g_\varphi=|\eta|^2+\left[dt^2+\left|-\frac12x\times dx-\frac12(b-x\xi)\right|^2+\frac{3}{4}|b|^2\right].$$
The metric on $B$ is obtained then by removing the vertical term $\eta$ so
\[
 g_{B,0}=dt^2+
 \left|-\frac12x\times dx-\frac12(b-x\xi)\right|^2
 +\frac34|b|^2.
\]
Define $$\alpha=\iota_V\omega_N,\qquad \beta=\iota_V\Psi_{N,+},\qquad \gamma=\iota_V\Psi_{N,-},\qquad \tau=\omega_N-\eta\wedge \alpha.$$
Then, we get that
$$\alpha=\frac{\sqrt3}{2}\xi,\qquad \beta=-\frac{\sqrt3}{4}d\xi,\qquad  \tau=-\frac{\sqrt3}{4}C,\qquad
 \gamma=\frac14S-\frac14E-\frac12U,$$
where
$$
 S=\frac12x\cdot(dx\times dx),\qquad
 E=\sum_{i=1}^3 dx_i\wedge b_i,\qquad U=\sum_{\mathrm{cyc}}x_i\, b_j\wedge b_k,\qquad
 C=\sum_{i=1}^3(x\times dx)_i\,\wedge b_i.
$$
Here the cross product of vector-valued $1$-forms uses exterior multiplication. The $2$-form $S$ is the unit-sphere area form with $\int_{S^2}S=4\pi$. 
The induced $SU(3)$ completion is
$$
 F_0=-dt\wedge\alpha+\beta,\qquad
 \psi_{+,0}=dt\wedge\tau-\alpha\wedge\gamma,\qquad
 \psi_{-,0}=-\alpha\wedge\tau-dt\wedge\gamma.\label{gtwo-explicit-base}
$$
Moreover, $$\kappa_0=d\eta=\frac12S-\frac12d\xi.$$
These forms verify \eqref{gtwo-lift1}--\eqref{gtwo-lift3} directly. Theorem~\ref{gtwo-recognition} gives pluricanonicity and non-integrability. Moreover, substitution in the homogeneous nearly K\"ahler structure gives
\begin{align*}
 \omega_N&=\eta\wedge\alpha+\tau,\\
 \Psi_{N,+}&=\eta\wedge\beta-\alpha\wedge\gamma,\\
 \Psi_{N,-}&=\eta\wedge\gamma+\alpha\wedge\beta.
\end{align*}
Finally, we remark that $\left[\frac{\kappa_0}{2\pi}\right]
 =\left[\frac{S}{4\pi}\right]\in H^2(B;\mathbb Z)$, because $\int_{S^2}S=4\pi$.

 \section{Pluricanonical metrics and He's flows}\label{flow-he}

In \cite{flow-He2021}, He considers the
evolution equations
\begin{equation}\label{flow-he-equations}
 \frac{\partial F}{\partial t}=-\Delta^{g(t)}F,
 \qquad
 \frac{\partial F}{\partial t}=-\delta^{g(t)}dF.
\end{equation}
For a fixed almost complex structure $J$, the evolving metric is determined
by the $J$-invariant part of $F$. Our solutions remain compatible with $J$,
so that $g(t)=F(t)(\cdot,J\cdot)$ throughout. We point out that these are
different equations: the eigenform identity for $\Delta^g$ does not imply
the same identity for $\delta^g d$.

We recall the following identities for a pluricanonical LCAK structure proved in \cite{addison2026pluricanonicallocallyconformallykahler}

\begin{equation}\label{flow-split-laplacian}
 d\delta^gF=(n-1)(\|\theta\|^2F-\theta\wedge J\theta),\qquad
 \delta^g dF=(n-1)\theta\wedge J\theta.
\end{equation}

A consequence of the pluricanonical LCAK identities is the following explicit solution of He's full Hodge equation.

\begin{thm}\label{flow-homothety}
Let $(M,g_0,J,F_0,\theta)$ be a pluricanonical LCAK manifold of real
dimension $2n$, $n\geq2$. Put $\lambda=(n-1)\|\theta\|_{g_0}^2$.
Then
\begin{equation}\label{flow-homothetic-solution}
 F(t)=(1-\lambda t)F_0,\qquad
 g(t)=(1-\lambda t)g_0
\end{equation}
is a solution of the full Hodge equation in
\eqref{flow-he-equations} for $1-\lambda t>0$. The almost complex
structure and the Lee form remain fixed, and the evolving structure is
pluricanonical. If $\lambda>0$, this family degenerates at
$t=\lambda^{-1}$; if $\lambda=0$, it is stationary.
\end{thm}

\begin{proof}
For a spatially constant $c(t)>0$ (it depends only on $t$), the Hodge Laplacian satisfies
$\Delta^{cg_0}=c^{-1}\Delta^{g_0}$ on differential forms. Consequently,
\[
 \Delta^{cg_0}(cF_0)=\Delta^{g_0}F_0=\lambda F_0.
\]
The equation for $F(t)=c(t)F_0$ is therefore $c'=-\lambda$.
Moreover, $dF(t)=\theta\wedge F(t)$, and constant metric homotheties
leave the Levi-Civita connection unchanged. Thus $(D\theta)^{J,+}=0$
remains true. The image of the Nijenhuis tensor remains $g$-orthogonal to
$T$ and $JT$, since multiplying the metric by a positive constant
preserves orthogonality.
\end{proof}

\subsection{The \texorpdfstring{$d^*d$}{d*d} equation}

We suppose that the manifold is connected and that $\theta$ is nowhere zero. Let
$ V=\operatorname{Span}\{T,JT\}$, where $T=\theta^{\sharp_{g_0}}$,
and let $ H=\mathcal V^{\perp_{g_0}}$. We write
\begin{equation}\label{flow-horizontal-vertical}
 \begin{split}
  F_V&=\|\theta\|^{-2}\theta\wedge J\theta,\qquad F_H=F_0-F_V,\\
  g_V&=\|\theta\|^{-2}(\theta\otimes\theta+J\theta\otimes J\theta),
  \qquad g_H=g_0-g_V.
 \end{split}
\end{equation}
The tensors $g_H$ and $g_V$ are the restrictions of $g_0$ to the
orthogonal distributions, extended by zero. The $2$-forms $F_H$ and
$F_V$ are their fundamental forms. Since $\|\theta\|$ is constant,
\begin{equation}\label{flow-block-differentials}
 F_H=-\|\theta\|^{-2}dJ\theta,\qquad dF_H=0,\qquad
 dF_V=\theta\wedge F_H.
\end{equation}

\begin{prop}\label{flow-canonical-deformation}
For any spatially constant $a(t),b(t)>0$, set
\begin{equation}\label{flow-two-parameter-family}
 F_{a,b}=aF_H+bF_V,\qquad g_{a,b}=ag_H+bg_V,
 \qquad\theta_{a,b}=\frac ba\theta.
\end{equation}
Then $(g_{a,b},J,F_{a,b},\theta_{a,b})$ is pluricanonical LCAK.
It satisfies
\begin{equation}\label{flow-family-identities}
 \|\theta_{a,b}\|_{g_{a,b}}^2=\frac b{a^2}\|\theta\|^2,
 \qquad
 \delta^{g_{a,b}}dF_{a,b}=\lambda\frac{b^2}{a^2}F_V,
 \qquad\lambda=(n-1)\|\theta\|^2.
\end{equation}
\end{prop}
\begin{proof}
Both blocks of $g_{a,b}$ are $J$-invariant, and $a,b>0$ give a positive definite
metric. From \eqref{flow-block-differentials},
\[
 dF_{a,b}=b\,\theta\wedge F_H=\theta_{a,b}\wedge F_{a,b},
 \qquad d\theta_{a,b}=0.
\]
The dual of $\theta_{a,b}$ for $g_{a,b}$ is $a^{-1}T$, which proves
the Lee norm formula and shows that its Lee plane is still
$ V$. Since $J$ is fixed, its Nijenhuis tensor is unchanged.
Its image lies in $ H$ for the initial pluricanonical structure,
and thus remains $g_{a,b}$-orthogonal to $ V$.

The remaining pluricanonical condition follows from
$$ dJ\theta_{a,b}=-\frac{b\|\theta\|^2}{a}F_H =-\|\theta_{a,b}\|_{g_{a,b}}^2F_{a,b} +\theta_{a,b}\wedge J\theta_{a,b}.$$

Indeed, the general LCAK identity
\[
 dJ\vartheta=-\|\vartheta\|^2F+\vartheta\wedge J\vartheta
       +2(D\vartheta)^{J,+}_{J\cdot,\cdot}+2N_{J\vartheta^\sharp},
\]
where $N_{X}(Y,Z)=g(N(Y,Z),X)$, from \cite[Proposition 2]{addison2026pluricanonicallocallyconformallykahler},
applied to the new metric, gives
$(D^{g_{a,b}}\theta_{a,b})^{J,+}=0$. Finally,
\eqref{flow-split-laplacian} for this pluricanonical structure gives
\[
 \delta^{g_{a,b}}dF_{a,b}
  =(n-1)\theta_{a,b}\wedge J\theta_{a,b}
  =\lambda\frac{b^2}{a^2}F_V.
\]
\end{proof}

\begin{thm}\label{flow-immortal-ddstar}
With the notation above, put
$b(t)=(1+\lambda t)^{-1}$. Then
\begin{equation}\label{flow-immortal-solution}
 F(t)=F_H+b(t)F_V,\qquad
 g(t)=g_H+b(t)g_V,\qquad
 \theta(t)=b(t)\theta,\qquad J(t)=J
\end{equation}
is a pluricanonical solution of
$$\partial_tF=-\delta^{g(t)}dF,$$ on $(-\lambda^{-1},\infty)$.
In particular, the solution is immortal in the forward direction.
On a compact manifold, it is the unique smooth solution with this
initial form in He's fixed-$J$ tamed class, on every common forward interval of
existence. It also solves He's $d^*d$ equation for compatible pairs with
the compatibility correction determined by the $J$-anti-invariant part
of $-\delta^{g(t)}dF(t)$.
\end{thm}
\begin{proof}
Proposition~\ref{flow-canonical-deformation} proves pluricanonicity at
every positive value of $b$. The flow equation implies that $a'(t)=0$, with $a(0)=1$. Then, for $a=1$, the evolution equation reduces to
$b'=-\lambda \, b^2$, $b(0)=1$.
The velocity $\dot{F}(t)$ is $J$-invariant, so the indicated compatibility correction
vanishes. On compact manifolds, uniqueness follows from
\cite[Theorem 2.1]{flow-He2021}, applied successively from any time at which
two smooth solutions agree.
\end{proof}

\bibliographystyle{ieeetr}

\bibliography{Geometry}

\end{document}